\documentclass[12pt,a4paper, twoside]{article}
\usepackage[latin1]{inputenc} 
\usepackage[T1]{fontenc} 
\usepackage[french]{babel}
\usepackage[sc]{mathpazo}
\usepackage{amsmath}
\usepackage{amsfonts}
\usepackage{amssymb}
\usepackage{amsthm}
\usepackage{mathrsfs}
\usepackage{geometry}
\usepackage{hyperref}
\usepackage{accents}
\usepackage{etoolbox}
\usepackage{stmaryrd}
\usepackage{caption}
\usepackage{mathtools}
\usepackage[mathscr]{euscript}
\usepackage{pgf,tikz,pgfplots,xcolor}
\usepackage{fancyhdr}
\usepackage[all]{xy}
\usetikzlibrary{cd}
\usetikzlibrary{arrows}
\pgfplotsset{compat=1.15}
\author{Emeryck Marie}
\newtheoremstyle{emi}{\topsep}{\topsep}{\itshape}{0pt}{}{\hspace{1ex}---}{1ex}{\thmname{{\itshape#1}}\thmnumber{ {\bfseries(#2.)}}\thmnote{ (#3.)}}
\theoremstyle{emi}
\newtheorem{theoreme}{Theorème}[subsection]

\newtheorem{proposition}[theoreme]{Proposition}
\newtheorem{definition}[theoreme]{Définition}
\newtheorem{corollaire}[theoreme]{Corollaire}

\newtheoremstyle{emiparagraphe}{\topsep}{\topsep}{}{0pt}{}{\hspace{1ex}---}{1ex}{\thmnumber{{\bfseries#2.}}\thmnote{ (#3.)}}
\theoremstyle{emiparagraphe}
\newtheorem{paragraphe}[theoreme]{}

\newtheoremstyle{emiremarque}{\topsep}{\topsep}{}{0pt}{}{\hspace{1ex}---}{1ex}{\thmname{{\itshape#1}}\thmnote{ (#3.)}}
\theoremstyle{emiremarque}
\newtheorem*{remarque}{Remarque.}
\newtheorem*{preuve}{Preuve.}
\AtEndEnvironment{preuve}{\qed}

\addto\captionsfrench{%
  \renewcommand{\contentsname}%
    {Contenu.}%
}

\addto\captionsfrench{}

\title{Continuité des racines d'après Rabinoff et Berkovich.}

\begin{document}

\maketitle

\begin{abstract}
	The content of this paper is a generalization of a theorem by Joseph Rabinoff: if $\mathscr{P}$ is a finite family of pointed and rational polyhedra in $N_\mathbb{R}$ such that there exists a fan in $N_\mathbb{R}$ that contains all the recession cones of the polyhedra of $\mathscr{P}$, if $k$ is a complete non-archimedean field, if $S$ is a connected and regular $k$-analytic space and $Y$ is a closed $k$-analytic subset of $U_{\mathscr{P}} \times_k S$ which is relative complete intersection and contained in the relative interior of $U_{\mathscr{P}} \times_k S$ over $S$, then the quasifiniteness of $\pi : Y \to S$ implies its flatness and finiteness; moreover, all the finite fibers of $\pi$ have the same length. This namely gives a analytic justification to the concept of stable intersection used in the theory of tropical intersection.
\end{abstract}

\tableofcontents

\textsc{Classification MSC :} 13, 14, 51, 52.

\textsc{Mots-clés :} Espaces de Berkovich, géométrie tropicale, géométrie polyédrale.

\section{Introduction et notations.}

\begin{center}
	\textsc{Introduction :}
\end{center}

	Ce texte trouve sa source dans l'article \cite{Rab} de Joseph Rabinoff dans lequel il établit, en utilisant la théorie des espaces rigides, une formule pour compter le nombres de zéros communs (comptés avec multiplicité et dont les valuations sont prescrites) de $n$ séries formelles convergentes à $n$ variable à coefficients dans un corps non archimédien non trivialement valué en terme de données de nature polyédrale --- liées au polyèdre de Newton des séries formelles. Cet article utilise le langage de la géométrie rigide ainsi que la théorie des modèles formels de Raynaud, le présent article se propose de reformuler les énoncés et les preuves dans le langage de la théorie des espaces de Berkovich, et par la même occasion de généraliser le théorème \cite[Theorem 9.8.]{Rab} que Rabinoff appelle  \og continuité des racines \fg{}  au cas d'une base de dimension finie arbitraire non nécessairement affinoïde mais régulière ; ce théorème de continuité des racines est l'ingrédient crucial pour passer du cas des polytopes traité par Osserman et Payne \cite[Corollary 5.2.4.]{OssPay} au cas d'un polyèdre --- et donne au passage une justification au principe d'\emph{intersection stable} utilisé en théorie de l'intersection tropicale. Le niveau de généralité établi ici permet de répondre positivement à la conjecture \cite[Remarque 9.11.]{Rab} de Rabinoff. Je tiens à remercier Antoine Chambert-Loir, qui m'a fait connaître cet article de Rabinoff, mais aussi pour la liberté qu'il m'a laissée lors de l'écriture de mon mémoire de Master 2 dont ce texte est l'augmentation d'une partie. Je tiens également à remercier les deux rapporteur·ses pour leurs remarques.

\begin{center}
	\textsc{Notations et conventions :}
\end{center}

Dans ce texte, $k$ désignera un corps complet non archimédien dont la valuation sera notée val. On note $\Gamma_k:=\mathrm{val}(k^\times)$ son groupe de valeurs et $\Gamma$ le groupe de valeurs d'une clôture algébrique de $k$. On définit une valeur absolue non archimédienne sur $k$ associée à cette valuation par:
\begin{center}
	$\forall x \in k, |x|:=\exp(-\mathrm{val}(x))$
\end{center}
avec la convention selon laquelle $\exp(-\infty):=0$. Réciproquement, si l'on dispose d'une valeur absolue non archimédienne $|\cdot|$ sur $k$, on peut y définir une valuation $v$ par :
\begin{center}
	$\forall x \in k, v(x):=-\ln(|x|)$
\end{center}
avec la convention selon laquelle $\ln(0):=+ \infty$.

On note $k^\mathrm{o}$ l'anneau de valuation de $k$, $k^\mathrm{oo}$ son idéal maximal et $\kappa:=k^\mathrm{o}/k^{\mathrm{oo}}$ son corps résiduel. \\
On fixe également un entier naturel $n$, un espace vectoriel réel $N_\mathbb{R}$ de dimension $n$ et un réseau $N$ dans $N_\mathbb{R}$, de sorte que $N_\mathbb{R}=N \otimes_\mathbb{Z} \mathbb{R}$. On note $M:=N^*$ le dual de $N$ et $M_\mathbb{R}:=N_\mathbb{R}^*$ le dual de $N_\mathbb{R}$. \\
La topologie dont est muni $\overline{\mathbb{R}}:=\mathbb{R} \cup \lbrace - \infty \rbrace$ est la topologie usuelle sur $\mathbb{R}$ à laquelle on adjoint une base de voisinages de $- \infty$ donnée par les intervalles fermés de la forme $([-\infty, a])_{a \in \mathbb{R}}$. Si $P$ est un polyèdre dans $N_\mathbb{R}$ et si $F \subseteq P$, on note $F \prec P$ pour dire que $F$ est une face de $P$. On notera aussi $\mathrm{Vert}(P)$ l'ensemble des sommets de $P$.

\section{Tropicalisation, géométrie polyédrale et géométrie analytique.}

\subsection{Géométrie polyédrale.}

Dans cette sous-section, on rappelle les constructions de géométrie polyédrale qui seront abondamment utilisées dans la suite ; on introduit également des constructions moins classiques comme la compactification d'un espace vectoriel le long d'un cône introduite par Kajiwara \cite{Kaj} et développée par Payne \cite{Pay}.

\begin{definition}[Polyèdre]
	Un \emph{polyèdre} $P$ dans $N_\mathbb{R}$ est une partie de $N_\mathbb{R}$ de la forme :
	 \begin{center}
	 	$P:=\lbrace v \in N_\mathbb{R} \mid \forall i \in \llbracket 1,r \rrbracket, \langle u_i, v \rangle \leq  a_i \rbrace$
	 \end{center}
	 où $r$ est un entier naturel non nul, $u_i$ une forme linéaire sur $N_\mathbb{R}$ et $a_i$ un nombre réel. On dit que :
	 \begin{itemize}
	 	\item $P$ est un \emph{polytope} si $P$ est un polyèdre compact.
	 	\item $P$ est un polyèdre \emph{rationnel} si les $u_i$ sont des éléments de $M$, le réseau dual.
	 	\item $P$ est $\Gamma$-affine (où $\Gamma \subseteq \mathbb{R}$) si les $a_i$ sont des éléments de $\Gamma$.
	 \end{itemize}
\end{definition}

\begin{remarque}
	Lorsque l'on écrit un polyèdre défini par les inégalités $(\langle u_i, \cdot \rangle \leq a_i)_{1 \leq i \leq r}$ on suppose toujours que les $u_i$ sont linéairement indépendants --- sur $\mathbb{Z}$ si le polyèdre est entier et sur $\mathbb{R}$ sinon.
\end{remarque}

\begin{definition}[Cône de récession]
	Soit $P$ un polyèdre dans $N_\mathbb{R}$ défini par $(\langle u_i, \cdot \rangle \leq a_i)_{1 \leq i \leq r}$. Le \emph{cône de récession de $P$ dans $N_\mathbb{R}$} est défini par :
	\begin{center}
		$\mathrm{Recc}(P):=\lbrace v \in N_\mathbb{R} \mid \forall i \in \llbracket 1,r \rrbracket, \langle u_i, v \rangle \leq 0 \rbrace$.
	\end{center}
	Il s'agit d'un cône dans $N_\mathbb{R}$.
\end{definition}

\begin{remarque}
	Si le polyèdre est supposé rationnel, alors son cône de récession est un cône \emph{rationnel}, c'est-à-dire qu'il admet un système générateur formé d'éléments du réseau $M$.
\end{remarque}

\begin{definition}[Cône/polyèdre pointé]
	On dit qu'un cône est \emph{pointé} si $\lbrace 0 \rbrace$ est l'une de ses faces. Un polyèdre est dit \emph{pointé} si son cône de récession est pointé.
\end{definition}

\begin{definition}[Cône polaire]
	Si $\sigma$ est un cône dans $N_\mathbb{R}$, on définit son \emph{cône polaire} par :
	\begin{center}
		$\sigma^\vee:=\lbrace u \in M_\mathbb{R} \mid \forall v \in \sigma, \langle u,v \rangle \leq 0 \rbrace$.
	\end{center}
	C'est un cône dans $M_\mathbb{R}$ qui est rationnel si $\sigma$ l'est.
\end{definition}

\begin{remarque}
	Lorsque $\sigma$ est un cône pointé dans $N_\mathbb{R}$, $\sigma^\vee$ n'est contenu dans aucun sous-espace strict de $M_\mathbb{R}$. 
\end{remarque}

\begin{paragraphe}
	Soit $\sigma$ un cône dans $N_\mathbb{R}$. La \emph{compactification partielle} de $N_\mathbb{R}$ le long de $\sigma$ est définie par :
	\begin{center}
		$N_\mathbb{R}(\sigma):=\mathrm{Hom}_{\mathbb{R}_+}(\sigma^\vee, \overline{\mathbb{R}})$.
	\end{center}
	On ne considère ici que les morphismes de monoïdes qui sont invariants sous l'action naturelle de $\mathbb{R}_+$ sur $\sigma^\vee$. On fait de $N_\mathbb{R}(\sigma)$ un espace topologique en le munissant de la topologie produit. Si de plus, on suppose $\sigma$ pointé, alors $N_\mathbb{R}(\lbrace 0 \rbrace)=N_\mathbb{R}$ s'injecte continûment dans $N_\mathbb{R}(\sigma)$
\end{paragraphe}

\begin{remarque}
	On adopte ici la convention selon laquelle $0 \cdot - \infty := 0$ ; par conséquent, le morphisme constant égal à $-\infty$ ne définit pas un élément de $N_\mathbb{R}(\sigma)$ puisqu'il n'est pas $\mathbb{R}_+$-equivariant. 
\end{remarque}

\begin{paragraphe}
	Si $\sigma^\vee:=\mathrm{Cone}(u_1, \dots, u_r)$, alors on a une immersion fermée topologique --- c'est-à-dire une application continue et injective qui est un homéomorphisme sur son image que l'on suppose fermée --- donnée par :
	\begin{center}
		$\iota : v \in N_\mathbb{R}(\sigma) \longmapsto (\langle u_i, v \rangle)_{1 \leq i \leq  r} \in \overline{\mathbb{R}}^r$.
	\end{center}
	On peut donc identifier $N_\mathbb{R}(\sigma)$ à un fermé de $\overline{\mathbb{R}}^r$.
\end{paragraphe}

\begin{definition}[Compactification d'un polyèdre pointé]
	Si $P$ est un polyèdre pointé dans $N_\mathbb{R}$, alors la \emph{compactification} $\overline{P}$ de $P$ est définie comme l'adhérence de $P$ dans $N_\mathbb{R}(\mathrm{Recc}(P))$.
\end{definition}

\begin{proposition}
	Si $P$ est un polyèdre de $N_\mathbb{R}$, alors $\overline{P}$ est compact.
\end{proposition}

\begin{preuve}
	Si $P$ est défini par les inégalités $(\langle u_i, \cdot \rangle \leq a_i)_{1 \leq i \leq r}$, alors à travers l'immersion $\iota$, $\overline{P}$ s'identifie au fermé $\prod_{i=1}^r [-\infty, a_i]$ de $\overline{\mathbb{R}}^r$ qui est compact.
\end{preuve}

La compactification de $N_\mathbb{R}$ le long d'un cône $\sigma$ revient en fait à compactifier $N_\mathbb{R}$ le long des faces de $\sigma$ comme le montre la proposition suivante dont on pourra trouver l'énoncé dans \cite[Proposition 3.19.]{Rab} :

\begin{proposition}
	Soit $P$ un polyèdre dans $N_\mathbb{R}$ et $\sigma$ un cône dans $N_\mathbb{R}$.
	\begin{enumerate}
		\item Pour toute face $\tau$ de $\sigma$ et tout $v \in N_\mathbb{R}/\mathrm{Span}(\tau)$, on définit $\iota(v) \in N_\mathbb{R}(\sigma)$ par :
		\begin{center}
			$\iota(v) : u \in \sigma^\vee \mapsto \begin{cases}
			\langle u,v \rangle & \text{si $u \in \tau^\perp$} \\
			-\infty & \text{sinon}
			\end{cases}$
		\end{center}
		Ceci induit une application :
		\begin{center}
			$\displaystyle \iota : \coprod_{\tau \prec \sigma} N_\mathbb{R}/\mathrm{Span}(\tau) \to N_\mathbb{R}(\sigma)$
		\end{center}
		qui est un homéomorphisme lorsque l'on munit son domaine de la topologie de la somme disjointe --- c'est-à-dire que pour tout $\tau \prec \sigma$, l'application $\iota|_{N_\mathbb{R}/\mathrm{Span}(\tau)}$ est une immersion topologique.
		\item Si $P$ est défini par les inégalités $(\langle u_i, \cdot \rangle)_{1 \leq i \leq r}$, alors on a, à travers l'identification (topologique) de 1. :
		\begin{center}
			$\displaystyle \overline{P}=\coprod_{\tau \prec \mathrm{Recc}(P)} \pi_{\tau}(P)$ où $\pi_{\tau} : N_\mathbb{R} \to N_\mathbb{R}/\mathrm{Span}(\tau)$.
		\end{center}
	\end{enumerate}
\end{proposition}

\begin{remarque}
	De la proposition ci-dessus, on déduit que si $P$ est déjà compact, alors $\overline{P} = P$ puisque dans ce cas $\mathrm{Recc}(P)=\lbrace 0 \rbrace$ et donc $\overline{P}=\pi_{\lbrace 0 \rbrace}(P)=P$. Cela peut être également prouvé directement puisque $N_\mathbb{R}(\lbrace 0 \rbrace)=N_\mathbb{R}$.
\end{remarque}

\begin{definition}[Complexe polyédral, éventail]
	Un \emph{complexe polyédral} est une famille finie $\Pi:=(P_i)_{1 \leq i \leq r}$ de polyèdres de $N_\mathbb{R}$ telle que :
	\begin{enumerate}
		\item $\Pi$ est stable par intersection deux-à-deux non vides.
		\item Toute face d'un élément de $\Pi$ est encore un élément de $\Pi$.
	\end{enumerate}
	Un élément de $\Pi$ est appelé une \emph{cellule} de $\Pi$ et le \emph{support} de $\Pi$ est défini comme l'union de ses cellules, on le note $|\Pi|$. Un complexe polyédral dont toutes les cellules sont des cônes est appelé un \emph{éventail} ; on dit qu'il est \emph{pointé} si toutes ses cellules le sont et on dit qu'il est \emph{complet} si son support est égal à $N_\mathbb{R}$.
\end{definition}

\begin{paragraphe}
	Si l'on dispose d'un éventail pointé $\Delta$ dans $N_\mathbb{R}$, on peut construire la compactification partielle de $N_\mathbb{R}$ le long de $\Delta$ en recollant les compactifications partielles $N_\mathbb{R}(\sigma)$ pour $\sigma$ cellule de $\Delta$ le long des immersions $N_\mathbb{R}(\tau) \hookrightarrow N_\mathbb{R}(\sigma)$ pour $\tau \prec \sigma$.
\end{paragraphe}

\begin{proposition}
	Si $\Delta$ est un éventail pointé dans $N_\mathbb{R}$, alors les assertions suivantes sont équivalentes :
	\begin{enumerate}
		\item $\Delta$ est un éventail complet.
		\item $N_\mathbb{R}(\Delta)$ est un espace topologique compact.
	\end{enumerate}
\end{proposition}

\subsection{L'application de tropicalisation.}

\begin{definition}
	Si $\sigma$ est un cône rationnel dans $N_\mathbb{R}$ et si l'on pose $S_\sigma:=\sigma^\vee \cap M$, alors la \emph{variété torique associée à $\sigma$} est définie comme le $k$-schéma:
	\begin{center}
		$X(\sigma):=\mathrm{Spec}(k[S_\sigma])$.
	\end{center}
\end{definition}

\begin{remarque}
		Puisque le cône est rationnel,  le lemme de Gordan (\cite[Section 1.2., Proposition 1]{Ful93}) implique que $X(\sigma)$ est un $k$-schéma de type fini.
\end{remarque}

On utilise le cône polaire pour définir $S_\sigma$ (là où la géométrie torique classique utilise le cône dual --- cf. \cite[Section 1.2., Proposition 1]{Ful93}) en raison de la définition de l'application de tropicalisation : on eût très bien pu la définir avec un signe moins, remplacer $- \infty$ par $+ \infty$ et utiliser le cône dual à la place du cône polaire ; on a simplement choisi d'utiliser la convention utilisée par Rabinoff.

\begin{paragraphe}
	Dans tout ce qui suit, si $X$ est un schéma de type fini sur un corps non archimédien $k$, $X^{an}$ désignera toujours son analytification au sens de Berkovich : il s'agit d'un espace $k$-analytique au sens de Berkovich, qui est sans bord. Si $X=\mathrm{Spec}(A)$ où $A$ est une $k$-algèbre de type fini, alors $X^{an}$ est l'ensemble des semi-normes multiplicatives sur $A$ étendant la valeur absolue de $k$. Plus de détails peuvent être trouvés dans \cite[§ 3.4. et 3.5.]{Ber1}.
\end{paragraphe}

On prouve à présent un résultat bien connu

\begin{definition}[Tropicalisation]
	Soit $\sigma$ un cône rationnel dans $N_\mathbb{R}$. \\
	On définit l'application de \emph{tropicalisation} par :
	\begin{center}
		$\mathrm{trop} : p \in X(\sigma)^{an} \mapsto \big( u \mapsto \ln(|x^u(p)|) \big) \in N_\mathbb{R}(\sigma)$.
	\end{center}
\end{definition}

\begin{proposition}
	La tropicalisation est une application continue, surjective et propre.
\end{proposition}

\begin{preuve}
	La continuité découle du fait que la fonction ln est continue sur $\mathbb{R}_+^\times$ et de la définition de la topologie sur $X(\sigma)^{an}$. La surjectivité découle du fait que l'on puisse construire une section à la tropicalisation, elle est donnée par $u \in N_\mathbb{R}(\sigma) \mapsto \eta_{e^{\langle u, \cdot \rangle}} \in X(\sigma)^{an}$ où $\eta_x$ désigne le point de Gauss en $x \in (\mathbb{R}^+)^n$ --- c'est-à-dire la semi-norme multiplicative définie par:
	\begin{center}
		$\displaystyle \eta_x : \sum_{u \in S_\sigma} a_u t^u \in k[S_\sigma] \mapsto \max\limits_{u \in S_\sigma}(|a_u|x^u) \in \mathbb{R}_+$.
	\end{center} 
	Concernant la propreté, puisque l'espace d'arrivée et celui de départ de l'application de tropicalisation sont localement compacts (le premier étant un bon espace $k$-analytique et le second homéomorphe à un fermé de $\overline{\mathbb{R}}^r$), il suffit de montrer que pour tout compact $C$ de $N_\mathbb{R}(\sigma)$, $\mathrm{trop}^{-1}(C)$ est compact. \\
	À présent, si $C$ est un compact de $N_\mathbb{R}(\sigma)$, alors il suffit de montrer que les domaines suivants sont compacts pour tout $s > 0$:
	\begin{center}
		$U_s:=\lbrace p \in X(\sigma)^{an} \mid \forall i \in \llbracket 1,r \rrbracket, |x_i(p)| \leq s \rbrace$.
	\end{center}
	En effet, par continuité de la tropicalisation, $\mathrm{trop}^{-1}(C)$ est un fermé de $X(\sigma)^{an}$ et est en fait borné puisque $C$ est borné dans $N_\mathbb{R}(\sigma)$ ; par conséquent, il existe $s>0$ tel que $\mathrm{trop}^{-1}(C)$ soit un fermé de $U_s$, qui est donc compact car $U_s$ l'est. Montrons à présent la compacité des $U_s$. \\
	Si $f:=\sum_{u \in S_\sigma} a_ux^u$ et que $p \in U_s$ pour $s>0$, alors :
	\begin{center}
		$|f(p)| \leq \max\limits_{u \in S_\sigma}(|a_u|\cdot |x^u(p)|) \leq |f(\eta_s)|$.
	\end{center}
	On en déduit alors qu'à travers l'immersion fermée topologique qui identifie une semi-norme $p \in X(\sigma)^{an}$ à l'ensemble de ses valeurs, $U_s$ est un fermé de $\prod_{f \in k[S_\sigma]}[0, |f(\eta_s)|]$ qui est compact en vertu du théorème de Tikhonov ; ainsi, $U_s$ est compact comme sous-espace fermé d'un comapct.
\end{preuve}

\subsection{Géométrie analytique non archimédienne à la Berkovich.}

\begin{definition}
	Si $P$ est un polyèdre rationnel pointé dans $N_\mathbb{R}$, le \emph{sous-espace polyédral associé à $P$} est défini par :
	\begin{center}
		$U_P:=\mathrm{trop}^{-1}(\overline{P})$.
	\end{center}
\end{definition}

\begin{remarque}
	Par propreté de la tropicalisation et compacité de $\overline{P}$, on déduit que $U_P$ est un sous-espace compact de $X(\mathrm{Recc}(P))^{an}$.
\end{remarque}

\begin{proposition}
	Si $P$ est un polyèdre rationnel pointé dans $N_\mathbb{R}$, alors $U_P$ est un domaine affinoïde de $X(\mathrm{Recc}(P))^{an}$.
\end{proposition}

\begin{preuve}
	Dans un premier temps, supposons que $P$ est un polytope ; en particulier $P=\overline{P}$. \\
	Écrivons alors $P:=\lbrace v \in N_\mathbb{R} \mid \forall i \in \llbracket 1,r \rrbracket, \langle u_i, v \rangle \leq a_i \rbrace$ où $(a_i)_{1 \leq i \leq r} \in \mathbb{R}_+^r$, $(u_i)_{1 \leq i \leq r} \in M^r$ et $\sigma:=\mathrm{Recc}(P)$, ainsi :
	\begin{center}
		$U_P=\lbrace p \in X(\sigma)^{an} \mid \forall i \in \llbracket 1,r \rrbracket, \sum_{j=1}^n u_{i,j} \ln(|x_j(p)|) \leq a_i \rbrace$ \\
		$=\lbrace p \in X(\sigma)^{an} \mid \forall i \in \llbracket 1,r \rrbracket, p(\prod_{j=1}^n x_j^{u_{i,j}}) \leq e^{a_i} \rbrace$.
	\end{center}
	où $(x_j)_j$ est un système générateur de la $k$-algèbre $k[S_\sigma]$ et $(u_{i,j})_j$ sont les coordonnées de $u_i$ pour tout $1 \leq i \leq r$ dans la base duale d'une base de $N_\mathbb{R}$. Ainsi, $U_P$ est un domaine de Weierstrass donc en particulier, un domaine affinoïde dans $X(\sigma)^{an}$. Pour le cas général, prendre l'adhérence dans $N_\mathbb{R}(\sigma)$ revient à considérer des points à l'infini, c'est-à-dire à autoriser la semi-norme à s'annuler sur les $x_j$ --- l'expression au-dessus reste valable.
\end{preuve}

\begin{paragraphe}
	Si $P$ est un polyèdre rationnel pointé dans $N_\mathbb{R}$ dont on note $\sigma$ le cône de récession, on considère :
	\begin{center}
		$\displaystyle k \langle U_P \rangle := \bigg\{ \sum_{u \in S_\sigma} a_ux^u \in k \llbracket S_\sigma \rrbracket \mid \forall v \in P, \lim_{|u| \to + \infty}(|a_u|e^{\langle u,v \rangle})=0 \bigg\}$.
	\end{center}
	C'est une $k$-algèbre de Banach pour la norme définie par $|f|_{\mathrm{sup}}:=\sup_{\xi \in U_P}(|f(\xi)|)$. \\
	On considère également la norme suivante $|f|_P:=\max_{(u,v) \in S_\sigma \times \mathrm{Vert}(P)}(|a_u|e^{\langle u,v \rangle})$.
\end{paragraphe}

\begin{theoreme}
	Soit $P$ un polyèdre rationnel pointé de $N_\mathbb{R}$ et $\sigma:=\mathrm{Recc}(P)$.
	\begin{enumerate}
		\item $k \langle U_P \rangle$ est une $k$-algèbre affinoïde pour la norme $|\cdot|_P$.
		\item Les deux normes définies ci-dessus coïncident.
		\item $k \langle U_P \rangle$ est une $k$-algèbre de Cohen-Macaulay.
	\end{enumerate}
\end{theoreme}

\begin{preuve}
	 Pour les deux premiers points, la preuve de \cite[Proposition 6.9.]{Rab} tient toujours à ceci près qu'il faut se ramener au cadre \emph{strictement} affinoïde ; pour cela, on utilise le procédé classique qui consiste à prendre le produit tensoriel complété avec l'algèbre $K_r$ pour $r$ un polyrayon convenable ; ceci ne change rien puisque ce foncteur est fidèlement exact et isométrique.
	
	Pour la troisième assertion, en vertu du théorème de Hochster (cf. par exemple \cite{Dan}), le schéma $X(\sigma)$ est de Cohen-Macaulay de dimension $n$ ainsi par \cite[Prop. 3.4.3.]{Ber1}, l'espace $k$-analytique $X(\sigma)^{an}$ est de Cohen-Macaulay de dimension $n$. Comme $U_P$ est un domaine affino\"{i}de de $X(\sigma)^{an}$, on déduit de \cite[Th. 3.4. B.]{Duc2} que $U_P$ est de Cohen Macaulay et donc que $k \langle U_P \rangle$ est une $k$-algèbre de Cohen-Macaulay.
\end{preuve}

\begin{remarque}
	La preuve ci-dessus montre que si $P$ est entier et $\Gamma$-affine, alors $k \langle U_P \rangle$ (resp. $U_P$) est une $k$-algèbre (resp. un domaine) \emph{strictement} affinoïde. Par ailleurs, lorsque $P$ est $\Gamma$-affine, on sait que son cône polaire $\sigma^\vee$ est de dimension maximale dans $M_\mathbb{R}$ : on en déduit donc que dans ce cas, $k \langle U_P \rangle$ est de dimension $n$.  
\end{remarque}

\section{Le théorème de continuité des racines.}

\subsection{Propreté en géométrie rigide et en théorie de Berkovich.}

L'article de Rabinoff utilise le langage de la géométrie rigide et des modèles formels de Raynaud. Dans cette section, la propreté des morphismes est une notion centrale, il s'agit alors dans un premier temps de comparer les définitions de la propreté dans le cadre de la géométrie rigide et dans le cadre de la théorie de Berkovich. 

\begin{definition}[Morphisme propre en géométrie rigide]
	Si $f : X \to Y$ est un morphisme entre deux espaces rigides, on dit que $f$ est \emph{propre} s'il est topologiquement séparé et qu'il existe un recouvrement admissible $(U_i)_{1 \leq i \leq s}$ (resp. $(V_i)_{1 \leq i \leq s}$) de $X$ (resp. $Y$) tels que $U_i$ soit relativement compact dans $V_i$ au-dessus de $Y$ pour tout $1 \leq i \leq s$.
\end{definition}

\begin{remarque} 
Cela signifie que pour tout $1 \leq i \leq s$, il existe une immersion fermée $V_i \hookrightarrow \mathbb{B}_k^n \times_k Y$ identifiant $U_i$ à $\mathbb{B}^n(r) \times_k Y$ pour $r \in ]0, 1[ \cap \Gamma$. Dans le cas d'espaces affinoïdes, il existe donc un épimorphisme (admissible) $\varphi : \Gamma(Y, \mathscr{O}_Y)\lbrace T_1, \dots, T_n \rbrace \twoheadrightarrow \Gamma(V_i, \mathscr{O}_{V_i})$ tel que $\rho(j(\varphi(T_m)))=r_m<1$ pour tout $1 \leq i \leq m$ où $j : \Gamma(X, \mathscr{O}_X) \to \Gamma(U, \mathscr{O}_U)$ est la restriction. En théorie de Berkovich, cela signifie exactement que le morphisme $j$ est intérieur pour $\Gamma(Y, \mathscr{O}_Y)$ : on peut donc voir la notion d'intérieur relatif comme la généralisation de la compacité relative introduite en géométrie rigide.
\end{remarque}

\begin{definition}[Morphisme propre]
	On dit qu'un morphisme $f : X \to Y$ entre deux espaces $k$-analytiques au sens de Berkovich est \emph{propre} s'il est topologiquement séparé, topologiquement propre et sans bord.
\end{definition}

En fait, ces deux définitions de propreté sont très liées:

\begin{proposition}[Temkin]
	Un morphisme entre espaces $k$-analytiques (au sens de Berkovich) séparés $f : X \to Y$ est propre au sens de Berkovich si et seulement si $f_0 : X_0 \to Y_0$ est propre au sens de la géométrie rigide. Ici, $X_0$ désigne l'ensemble des \emph{points rigides} de $X$, c'est-à-dire les points $x \in X$ tels que $[\mathscr{H}(x):k]$ est fini.
\end{proposition}

\begin{preuve}
	Le résultat et sa preuve se trouvent dans \cite[Corollary 4.5.]{Tem3}.
\end{preuve}

\subsection{Le théorème de continuité des racines.}

\begin{theoreme}[Théorème de Kiehl]
	Si $f : X \to Y$ est un morphisme propre entre deux espaces $k$-analytiques et que $\mathscr{F}$ est un $\mathscr{O}_X$-module cohérent, alors $f_*\mathscr{F}$ est un $\mathscr{O}_Y$-module cohérent.
\end{theoreme}

\begin{corollaire}
	Si $f : \mathscr{M}(\mathcal{B}) \to \mathscr{M}(\mathcal{A})$ est un morphisme propre, alors $\mathcal{B}$ est une $\mathcal{A}$-algèbre finie.
\end{corollaire}

\begin{proposition}
	Soit $P' \subseteq P$ deux polyèdres rationnels pointés dans $N_\mathbb{R}$ tels que $\mathrm{Recc}(P')$ soit une face du cône $\mathrm{Recc}(P)$. Si $P' \subseteq \mathrm{Relint}(P)$, alors $U_{P'} \subseteq \mathrm{Int}(U_P)$.
\end{proposition}

\begin{preuve}
	La preuve est similaire à \cite[Lemma 9.5.]{Rab}, en changeant la terminologie de la géométrie rigide par celle de la théorie de Berkovich, c'est-à-dire la compacité relative par des intérieurs relatifs.	 
\end{preuve}

\begin{proposition}[Critère tropical de finitude]
	Soit $P$ un polyèdre pointé dans $N_\mathbb{R}$. \\
	Si $Y$ est un fermé analytique de $U_P$ dont la tropicalisation est contenue dans $\mathrm{Relint}(\overline{P})$, alors $Y \to \mathscr{M}(k)$ est fini.
\end{proposition}

\begin{preuve}
	Soit $I$ l'idéal de $k \langle U_P \rangle$ définissant $Y$ comme fermé de $U_P$. En écrivant $P:=\lbrace v \in N_\mathbb{R} \mid \forall i \in \llbracket 1,r \rrbracket, \langle u_i, v \rangle \leq a_i \rbrace$, l'hypothèse sur $Y$ implique qu'il existe $b_i<a_i$ tel que $\mathrm{trop}(Y)$ est contenu dans l'adhérence du polyèdre $P'$ défini par les inégalités $(\langle u_i, \cdot \rangle \leq b_i)_{1 \leq i \leq r}$ et l'on peut prendre $P'$ de même cône de récession que $P$. La proposition 3.2.3. implique que $U_{P'} \subseteq \mathrm{Int}(U_P)$ et l'on a alors :
	\begin{center}
		$\mathrm{Int}(Y)=\mathrm{Int}(Y/U_P) \cap ( \mathrm{Int}(U_P) \cap Y )=Y \cap \mathrm{Int}(U_P) \supseteq Y$.
	\end{center}
	La première égalité vient du troisième point de \cite[Proposition 2.5.8.]{Ber1}, la seconde vient du fait que $Y \hookrightarrow U_P$ est une immersion fermée et donc un morphisme fini, ce qui implique que $\mathrm{Int}(Y/U_P)=Y$. L'inclusion finale vient du fait que  $Y \subseteq U_{P'} \subseteq \mathrm{Int}(U_P)$. On en déduit alors que $\partial Y = \emptyset$ donc $Y$ est un espace $k$-affinoide propre et $Y \to \mathscr{M}(k)$ est fini en vertu du corollaire 3.2.2.
\end{preuve}

\begin{paragraphe}
	Si l'on dispose d'une famille finie $(P_i)_{1 \leq i \leq r}$ de polyèdres pointés dans $N_\mathbb{R}$, alors leurs cônes de récession ont $0$ comme sommet commun. On aimerait que leurs cônes de récession respectifs soient des cônes d'un même éventail $\Delta$. Ce n'est malheureusement pas toujours possible : si l'on prend deux cônes pointés qui s'intersectent, l'éventail recherché devrait comporter trois cônes mais il n'y a que deux polyèdres. On introduit donc la définition suivante:
\end{paragraphe}

\begin{definition}[Polyèdres simultanément compactifiables]
	On dit qu'une famille $(P_i)_{i \in I}$ de polyèdres dans $N_\mathbb{R}$ est \emph{simultanément compactifiable} s'il existe un éventail de $N_\mathbb{R}$ dont les cônes de récession des $P_i$ soient des faces. 
\end{definition}

\begin{remarque}
	C'est par exemple le cas s'il existe $i_0 \in I$ tel que pour tout $i \in I$, $\mathrm{Recc}(P_i)$ est une face de $\mathrm{Recc}(P_{i_0})$.
\end{remarque}

\begin{paragraphe}
	Si $\mathscr{P}$ est une famille de polyèdres de $N_\mathbb{R}$ simultanément compactifiables dans un éventail noté $\Delta$, on pose $U_{\mathscr{P}}:=\mathrm{trop}^{-1}(\cup_{P \in \mathscr{P}} \overline{P})$, c'est un domaine analytique dans $X(\Delta)^{an}$ en tant qu'union \emph{finie} de domaines affinoïdes dans $X(\Delta)^{an}$.
\end{paragraphe}

\begin{theoreme}[Théorème de continuité des racines, version globale]
	Soit $S$ un espace $k$-analytique régulier et connexe. Soit $\mathscr{P}$ une famille finie de polyèdres pointés et rationnels dans $N_\mathbb{R}$ simultanément compactifiables dans un éventail noté $\Delta$. Soit $Y$ un fermé analytique de $U_\mathscr{P} \times_k S$ tel que $Y \subseteq \mathrm{Int}(U_\mathscr{P} \times_k S/S)$ tel que pour tout $s \in S$ et tout voisinage affino\"{i}de $V$ de $s$ dans $S$, le fermé $Y \cap U_P \times_k V$ soit défini comme le lieu des zéros de $d:=\dim_{\mathrm{Krull}}(k \langle U_P \rangle)$ éléments de $\Gamma(V, \mathscr{O}_V) \widehat{\otimes}_k k \langle U_P \rangle$ pour tout $P \in \mathscr{P}$. Alors :
	\begin{enumerate}
		\item Si $\pi : Y \to S$ est quasifini, alors $\pi$ est un morphisme fini et plat.
		\item Toutes les fibres finies de $\pi$ ont la même longueur.
	\end{enumerate}
\end{theoreme}

\begin{remarque}
	La seconde assertion est vraie y compris lorsque $\pi$ n'est pas quasifini.
\end{remarque}

\begin{preuve}
	Prouvons d'abord la propreté de $\pi$. On sait déjà que $\pi$ est topologiquement séparé puisque c'est la restriction d'une projection ; prouvons à présent que $\pi$ est topologiquement propre. Si $K$ est un compact de $S$, alors par compacité de $U_{\mathscr{P}}$, $U_{\mathscr{P}} \times_k K$ est compact ; par ailleurs, $Y$ est fermé dans $U_{\mathscr{P}} \times_k K$ donc $\pi^{-1}(K)$ est compact donc $\pi$ est topologiquement propre. Concernant le bord, étant donné que l'on a supposé $Y \subseteq \mathrm{Int}(U_{\mathscr{P}} \times_k S/S)$ et puisque $Y \hookrightarrow U_{\mathscr{P}} \times_k S$ est une immersion fermée, on a :
	\begin{center}
		$\mathrm{Int}(Y/S)=Y \cap \mathrm{Int}(U_{\mathscr{P}} \times_k S/S) = Y$.
	\end{center}
	Cela signifie précisément que $\partial(Y/S)=\emptyset$ et que $\pi$ est donc un morphisme propre. Puisque $\pi$ est également quasifini (par hypothèse), on en déduit que $\pi$ est un morphisme fini par \cite[Prop. 3.3.8]{Ber1}. \\
	Par définition (\cite[Definition 4.1.8.]{Duc3}), la platitude se vérifie sur un \emph{bon} domaine analytique : il suffit de prouver le résultat lorsque $S=\mathscr{M}(\mathcal{A})$ où $\mathcal{A}$ est une $k$-algèbre affinoïde régulière et de considérer $U_P \times_k S \cong \mathscr{M}(\mathcal{A} \widehat{\otimes}_k k \langle U_P \rangle)$ lorsque $P \in \mathscr{P}$. \\
	Écrivons $I:=(f_1, \dots, f_d)$ l'idéal définissant $Y \cap (U_P \times_k S)$ et $\mathcal{B}:=(\mathcal{A} \widehat{\otimes}_k k \langle U_P \rangle)/I$. \\
	En général, la platitude en théorie de Berkovich est différente de celle définie classiquement en géométrie algébrique \footnote{On impose en fait la stabilité par changement de base dans la définition --- cf. \cite[Ch. 4.]{Duc3}.} puisque l'on travaille avec des produits tensoriels complétés, la platitude au sens de la géométrie algébrique n'est pas stable par changement de base ; toutefois, pour un morphisme fini, les deux notions coïncident comme le montre \cite[Proposition 4.3.1.]{Duc3} --- c'est pourquoi la finitude de $\pi$ a été prouvée avant --- : il suffit alors de montrer que $\mathcal{B}$ est une $\mathcal{A}$-algèbre plate. \\
	Par \cite[Lemma 2.1.2.]{Ber2}, $\mathcal{A}\langle U_P \rangle:= \mathcal{A} \widehat{\otimes}_k k \langle U_P \rangle$ est une $\mathcal{A}$-algèbre plate. En observant qu'un anneau est de Cohen-Macaulay si et seulement si il vérifie la condition $S_n$ pour tout $n$. En appliquant \cite[Theorem 11.3.3.]{Duc1} aux faisceaux structuraux $\mathscr{O}_{\mathscr{M}(\mathcal{A})}$ et $\mathscr{O}_{\mathscr{M}(\mathcal{A}\langle U_P \rangle)}$, les fibres étant Cohen-Macaulay par le troisième point du théorème 2.3.4. et la platitude de la $\mathcal{A}$-algèbre $\mathcal{A} \langle U_P \rangle$ a été mentionnée auparavant, on déduit que $\mathcal{A} \langle U_P \rangle$ est également de Cohen-Macaulay. Par ailleurs, on sait que la dimension de Krull de $\mathcal{A} \langle U_P \rangle$ est égale à $d+\dim(\mathcal{A})$. Puisqu'un anneau de Cohen-Macaulay est caténaire, on a l'égalité suivante sur les dimensions de Krull :
	\begin{center}
		$\dim(\mathcal{B})=d+\dim(\mathcal{A})-\mathrm{ht}(I)$.
	\end{center}
	Par le \emph{Hauptidealsatz} de Krull, on a $\mathrm{ht}(I) \leq d$ donc $\dim(\mathcal{B}) \geq \dim(\mathcal{A})$ mais $\pi$ est quasifini donc $\dim(\mathcal{B}) \leq \dim(\mathcal{A})$ donc $\dim(\mathcal{B})=\dim(\mathcal{A})$ ; on en déduit alors que $\dim(\mathcal{A}\langle U_P \rangle)-\dim(\mathcal{B})=d$ et puisque $\mathcal{A} \langle U_P \rangle$ est de Cohen-Macaulay, on en déduit que $\mathcal{B}$ est également de Cohen-Macaulay. En appliquant le \textit{Miracle flatness} de Matsumura (\cite[Theorem 23.1.]{Mat}) sur les anneaux locaux, on en déduit que $\pi$ est un morphisme plat. 
	
	Concernant le second point, on peut toujours se ramener au cas où $\pi$ est un morphisme quasifini --- cas dans lequel où l'assertion découle de la platitude de $\pi$. En effet, en vertu du théorème de semi-continuité de la fibre \cite{Duc3}, l'ensemble des points de $S$ qui ne sont pas isolés dans leur fibre est un fermé $Z$ de $Y$. Par la preuve du premier point, $\pi$ est propre donc par le théorème de l'application de Remmert (\cite[Prop. 3.3.6.]{Ber1}), $\pi(Z)$ est un fermé $k$-analytique de $S$ et puisque $S$ est localement noethérien et normal, c'est une union disjointe d'espaces $k$-analytiques normaux irréductibles et puisque $S$ est connexe, on en déduit que $S$ est irréductible donc par \cite[Cor. 3.3.20.]{Ber1}, $S-\pi(Z)$ est encore connexe et par construction, $\pi|_{\pi^{-1}(S-\pi(Z))}$ est quasifini.
\end{preuve}

\begin{remarque}
	La condition $Y \subseteq \mathrm{Int}(U_\mathscr{P} \times_k S/S)$ est plus générale que celle demandée dans le théorème initial de Rabinoff \cite[Th. 9.8]{Rab} puisqu'avec ses notations, si l'on a $P_i' \subseteq \mathrm{Relint}(P_i)$ pour tout $1 \leq i \leq r$, alors on a $U_{P_i'} \subseteq \mathrm{Int}(U_{P_i})$ par la proposition 3.2.4. et par \cite[Lemma 4.1.4.]{Ber1}, on a :
	\begin{center}
		$U_{P_i'} \times_k S = \mathrm{pr}_1^{-1}(U_{P_i'}) \subseteq \mathrm{pr}_1^{-1}(\mathrm{Int}(U_{P_i})) \subseteq \mathrm{Int}(U_{P_i} \times_k S/S)$.
	\end{center}
	Puisque $Y$ est un fermé de $U_{\mathscr{P}'} \times_k S$ et que l'intérieur d'un morphisme peut être déterminé sur un recouvrement affinoïde fini, on en déduit effectivement que $Y \subseteq \mathrm{Int}(U_{\mathscr{P}} \times_k S/S)$. \\
	La condition sur le nombre d'équations du fermé $Y$ est également plus générale puisque si $P$ est $\Gamma$-affine, alors $k \langle U_P \rangle$ est une algèbre \emph{strictement} $k$-affino\"{i}de donc sa dimension de Krull est le nombre de formes linéaires définissant $P$ d'après le théorème 2.3.4. 
\end{remarque}

\begin{remarque}
\begin{itemize}
	\item Lorsque $\dim(S)=1$, le théorème 3.2.8. reste vrai lorsque $S$ est seulement connexe et \emph{normal} puisqu'en raison de l'hypothèse sur la dimension, il suffit de montrer que (en reprenant les notations de la preuve) $\mathcal{B}$ est une $\mathcal{A}$-algèbre sans torsion, ce qui résulte de l'\emph{unmixedness theorem}.
	\item Lorsque $P$ est $\Gamma$-affine, l'algèbre $k \langle U_P \rangle$ est \emph{strictement} $k$-affinoïde donc sa dimension de Krull coïncide avec sa dimension $k$-analytique --- au sens de \cite[§ 1.4.]{Duc3} ; on en déduit que dans ce cas, la dimension de l'espace $k$-analytique $U_P$. En général, la dimension de Krull de $k \langle U_P \rangle$ est seulement majorée par la dimension de l'espace $k$-analytique $U_P$, c'est-à-dire la dimension de $X(\mathrm{Recc}(P))^{an}$ par \cite[Lemme 1.15.]{Duc4} ou bien même la dimension de la variété torique affine $X(\mathrm{Recc}(P))$ --- par \cite[Theorem 3.4.8.]{Ber1} pour le cadre non trivialement valué et \cite[Theorem 3.5.3.]{Ber1} pour le cadre trivialement valué --- c'est-à-dire $n$ puisque $P$ est supposé pointé.
\end{itemize}
	
\end{remarque}

La preuve du théorème 3.2.8. met en exergue un résultat de nature locale qu'il est peut-être bon d'expliciter.

\begin{proposition}[Théorème de continuité des racines, version locale]
	Soit $\mathcal{A}$ une $k$-algèbre affinoïde régulière. On note $S:=\mathscr{M}(\mathcal{A})$ et $X:=\mathscr{M}(\mathcal{B})$ un espace $k$-affinoïde de Cohen-Macaulay et de dimension $d+\dim(S)$. Soit $(f_i)_{1 \leq i \leq d} \in \mathcal{B}^n$ et $Y:=V(f_1, \dots, f_d)^{an}$. S'il existe un morphisme $\pi : Y \to S$ tel qu'il existe $t \in S$ pour lequel $Y_t:=\pi^{-1}(t) \cap Y$ est de dimension nulle et n'étant pas l'image d'un point du bord de $S$, alors il existe un voisinage affinoïde $V$ de $t$ dans $S$ tel que $\pi|_{\pi^{-1}(V)}$ soit un morphisme fini et plat. En particulier, pour tout $s \in V$, $Y_s$ est fini et de même longueur que $Y_t$.
\end{proposition}

\begin{preuve}
	La preuve se ramène essentiellement à celle de 3.2.8., à ceci près que l'on doit ajouter l'hypothèse sur le fait que $t$ ne soit pas dans l'image d'un point du bord de $S$, cela assure la propreté de $\pi|_{\pi^{-1}(V)}$ et donc sa finitude ; une référence possible est \cite[Proposition 3.1.4.]{Ber2}.
\end{preuve}

\begin{remarque}
	De manière informelle, on peut penser à ce résultat comme une constance locale du nombre de zéros avec multiplicité et donc une \og continuité des racines \fg{} comme expliqué dans \cite[Example 10.3.]{Rab}. 
\end{remarque}

Du théorème global de continuité des racines et du critère tropical de finitude, on peut déduire un résultat plus simple à mettre en place que le premier cité :

\begin{corollaire}
	Soit $\mathscr{A}$ une algèbre $k$-affinoïde régulière, $P$ un polyèdre pointé et rationnel de $N_\mathbb{R}$, $(f_i)_{1 \leq i \leq d} \in \mathscr{A} \langle U_P \rangle^d$ avec $d:= \dim_{\mathrm{Krull}}(k \langle U_P \rangle)$. Si l'on considère $Y:=V(f_1, \dots, f_d) \subseteq U_P \times_k \mathscr{M}(\mathscr{A})$ et que $\mathrm{trop}(Y_s) \subseteq \mathrm{Relint}(\overline{P})$, alors le morphisme $Y \to \mathscr{M}(\mathscr{A})$ est fini et plat.
\end{corollaire}

\begin{remarque}
	Ici, on ne demande rien sur la finitude des fibres, cette dernière est impliqué par le critère tropical de finitude.
\end{remarque}

Détaillons à présent un exemple pour illustrer ce corollaire :

\begin{paragraphe}
	Prenons $n=2$, $p$ un nombre premier et $k:=\mathbb{Q}_p$ muni de la valuation $p$-adique notée $v_p$. Considérons le polyèdre rationnel pointé $P:=[-3,-1] \times ]- \infty, 0]$, son cône de récession est donné par $\sigma=\mathrm{Cone}((0,-1))$ qui est effectivement un cône pointé. On considère les deux équations suivantes :
	\begin{center}
		$f_1(x,y,t_1, t_2):=t_2+x+t_1y$ et $f_2(x,y,t_1,t_2)=p^2+x+y$.
	\end{center}
	Considérons à présent $Z:=V(f_1, f_2) \subseteq X(\sigma) \times \mathbb{A}^2_{\mathbb{Q}_p}$. Pour $(t_1, t_2)$ fixé, dessinons les tropicalisations des hypersurfaces des tropicalisations $V(f_{1, t_1, t_2})$ et $V(f_{1, t_1, t_2})$ dans $N_\mathbb{R}(\sigma)=\mathbb{R} \times [- \infty, + \infty[$ ; on obtient la figure suivante pour $(v_p(t_1), v_p(t_2))=(-8, 6)$
	\begin{center}
		\definecolor{ffqqqq}{rgb}{1,0,0}
\definecolor{qqffqq}{rgb}{0,1,0}
\definecolor{cqcqcq}{rgb}{0.75,0.75,0.75}
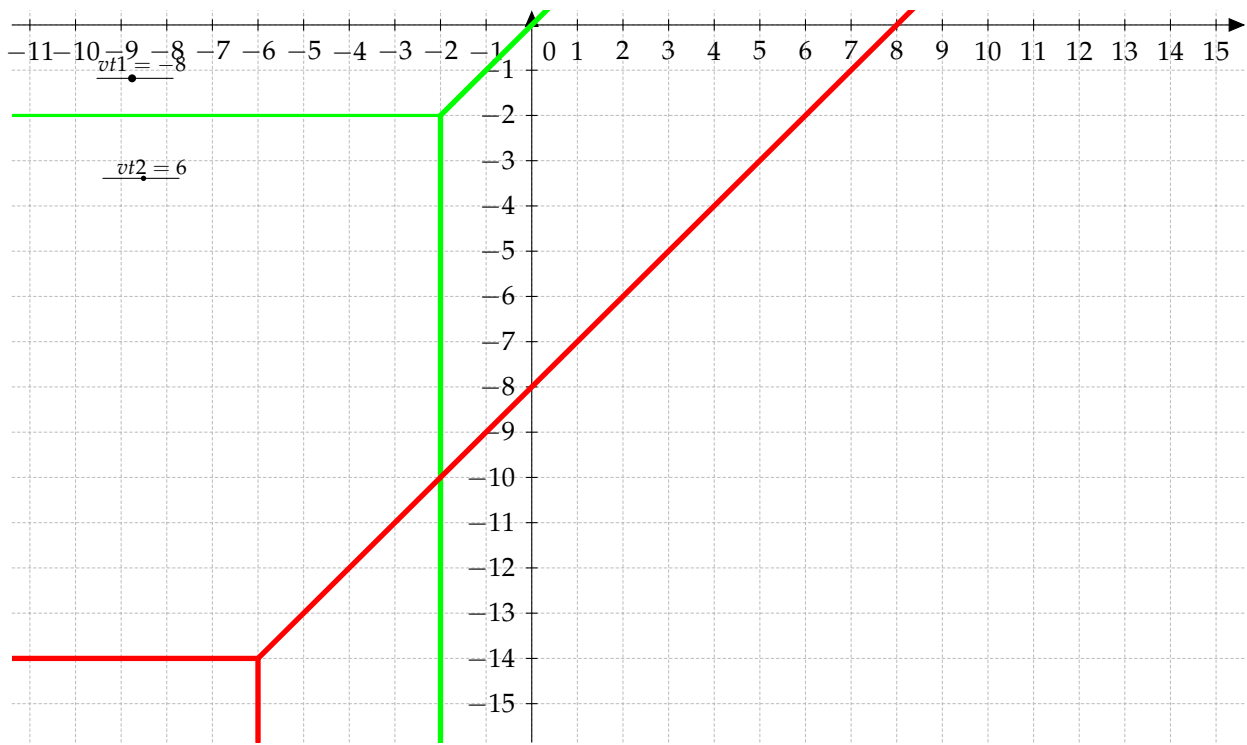
\begin{tikzpicture}[line cap=round,line join=round,>=triangle 45,x=0.6cm,y=0.6cm]
\draw [color=cqcqcq,dash pattern=on 1pt off 1pt, xstep=0.6cm,ystep=0.6cm] (-11.39,-15.85) grid (15.64,0.32);
\draw[->,color=black] (-11.39,0) -- (15.64,0);
\foreach \x in {-11,-10,-9,-8,-7,-6,-5,-4,-3,-2,-1,1,2,3,4,5,6,7,8,9,10,11,12,13,14,15}
\draw[shift={(\x,0)},color=black] (0pt,2pt) -- (0pt,-2pt) node[below] {\footnotesize $\x$};
\draw[->,color=black] (0,-15.85) -- (0,0.32);
\foreach \y in {-15,-14,-13,-12,-11,-10,-9,-8,-7,-6,-5,-4,-3,-2,-1}
\draw[shift={(0,\y)},color=black] (2pt,0pt) -- (-2pt,0pt) node[left] {\footnotesize $\y$};
\draw[color=black] (0pt,-10pt) node[right] {\footnotesize $0$};
\clip(-11.39,-15.85) rectangle (15.64,0.32);
\draw(-9.53,-1.18) -- (-7.86,-1.18);
\draw[line width=0.4pt] (-9.4,-3.39) -- (-7.73,-3.39);
\draw [line width=2pt,color=qqffqq,domain=-2.0:15.642847381912082] plot(\x,{(-0--202*\x)/202});
\draw [line width=1.2pt,color=qqffqq,domain=-11.392105288108946:-2.0] plot(\x,{(--188.25-0*\x)/-94.12});
\draw [line width=2pt,color=qqffqq] (-2,-2) -- (-2,-15.85);
\draw [line width=2pt,color=ffqqqq,domain=-6.0:15.642847381912082] plot(\x,{(-1156.41--144.55*\x)/144.55});
\draw [line width=2pt,color=ffqqqq] (-6,-14) -- (-6,-15.85);
\draw [line width=2pt,color=ffqqqq,domain=-11.392105288108946:-6.0] plot(\x,{(--3013.12-0*\x)/-215.22});
\begin{scriptsize}
\fill [color=black] (-8.76,-1.18) circle (1.5pt);
\draw[color=black] (-8.55,-0.88) node {$vt1 = -8$};
\fill [color=black] (-8.51,-3.39) circle (1.0pt);
\draw[color=black] (-8.33,-3.13) node {$vt2 = 6$};
\draw[color=qqffqq] (99.34,98.92) node {$r21$};
\draw[color=qqffqq] (-49.16,-1.87) node {$r22$};
\draw[color=qqffqq] (-1.83,-48.19) node {$r23$};
\draw[color=ffqqqq] (66.62,58.19) node {$r11$};
\draw[color=ffqqqq] (-5.83,-107.85) node {$r13$};
\draw[color=ffqqqq] (-113.71,-13.87) node {$r12$};
\end{scriptsize}
\end{tikzpicture}
	\end{center}
	où la tropicalisation de $V(f_{1, t_2, t_2})^{an}$ est dessinée en rouge et celle de $V(f_{2, t_2, t_2})^{an}$ est dessinée en vert. L'effet de la variation du couple $(v_p(t_1), v_p(t_2))$ sur les tropicalisations respectives des spécialisations est le suivant:
	\begin{itemize}
		\item Lorsque $v_p(t_1)$ augmente (resp. diminue), le graphe rouge est translaté vers le haut (resp. bas).
		\item Lorsque $v_p(t_2)$ augmente (resp. diminue), le graphe rouge est translatée selon le vecteur $(1,1)$ (resp. $(-1,-1)$).
\end{itemize}
On peut alors remarquer que si l'on ne fait pas trop varier la valuation des $t_i$ du couple $(-8, 6)$, la tropicalisation de la fibre $Z^{an}_{(t_1, t_2)}$  --- qui n'est autre que l'intersection des deux tropicalisations --- reste contenue dans l'intérieur de $\overline{P}=[-3,-1] \times [- \infty, 0] \subseteq N_\mathbb{R}(\sigma)$. On peut donc appliquer le corollaire 3.2.11. pour $\mathscr{A}$ convenable \footnote{typiquement, l'algèbre des fonctions analytiques sur une poly-couronne --- elle sera ici même \emph{strictement} $\mathbb{Q}_p$-affinoïde puisque les scalaires pourront être pris dans le groupe de valeurs de la valuation $v_p$ ; en particulier, elle sera bien de dimension de Krull égale à deux et son spectre satisfait les hypothèses du corollaire 3.2.10.} et $Y:= Z^{an} \cap (U_P \times \mathscr{M}(\mathscr{A}))$ et l'on obtient en particulier que toutes les fibres au-dessus de $\mathscr{M}(\mathscr{A})$ ont la même longueur, y compris celle pour $v_p(t_2)=2$ où les deux tropicalisation s'intersectent en une demi droite. 
\end{paragraphe}

\begin{remarque}
	De manière plus générale, on voit que ce résultat de continuité globale des racines permet de justifier le concept d'\emph{intersection stable} \footnote{consistant à faire varier un peu les deux courbes afin qu'elle ne se superposent plus --- cf. \cite[Theorem 4.3]{SRT03} par exemple.} qui est fondamental en théorie de l'intersection tropicale --- qui peut-être justifié en utilisant la condition d'équilibre sur les variétés tropicales. On pourra se reporter à la section 12 de \cite{Rab} et plus précisément \cite[Definition 12.7.]{Rab} pour la définition de l'intersection stable ainsi que \cite[Corollary 12.12.]{Rab} pour la justification de la légitimité de cette définition.
\end{remarque}

\bibliography{biblio}
\bibliographystyle{plain}

\textsc{Emeryck Marie}, Technische Universität Chemnitz, Fakultät für Mathematik ; Reichenhainer Stra\ss e 39, 09126 Chemnitz, Germany $\bullet$ \textsc{E-mail} : \texttt{\href{mailto:emeryck.marie@mathematik.tu-chemnitz.de}{emeryck.marie@mathematik.tu-chemnitz.de}}
\end{document}